\documentclass[english]{cccconf}
\usepackage[comma,numbers,square,sort&compress]{natbib}
\usepackage{amsmath,amssymb,amsfonts}
\allowdisplaybreaks[4]
\usepackage{epstopdf}
\usepackage{subfigure}
\newtheorem{thm}{\textbf{Theorem}}
\newtheorem{assu}{\textbf{Assumption}}
\newtheorem{rem}{\textbf{Remark}}
\newtheorem{lemma}{\textbf{Lemma}}
\newtheorem{corollary}{\textbf{Corollary}}

\begin{document}

\title{Online Distributed  Optimization with Clipped Stochastic Gradients: High Probability Bound of Regrets}

\author{Yuchen Yang\aref{amss},
        Kaihong Lu\aref{hit},
        Long Wang\aref{amss}}

\affiliation[amss]{Center for Systems and Control, College of
			Engineering,
        Peking University, Beijing 100871, P.~R.~China
        \email{ycyang@stu.pku.edu.cn, longwang@pku.edu.cn}}
\affiliation[hit]{College of Electrical Engineering and Automation, Shandong University of Science and Technology,  Qingdao 266590, P.~R.~China
        \email{khong$\_$lu@163.com}}

\maketitle

\begin{abstract}
In this paper, the problem of distributed optimization is studied via a network of agents. Each agent only has access to a stochastic gradient of its own objective function in the previous time, and can communicate with its neighbors via a network. To handle this problem, an online distributed clipped stochastic gradient descent algorithm is proposed. Dynamic regrets are
used to capture the performance of the algorithm. Particularly, the high probability bounds of regrets are analyzed when the
stochastic gradients satisfy the heavy-tailed noise condition. For the convex case, the offline benchmark of the dynamic regret is to seek the minimizer of the objective function each time. Under mild assumptions on the graph connectivity, we prove that the dynamic regret grows sublinearly with high probability under a certain clipping parameter. For the non-convex case, the offline benchmark of the dynamic regret is to find the stationary point of the objective function each time. We show that the dynamic regret grows sublinearly with high probability if the variation of the objective function grows within a certain rate. Finally, numerical simulations are provided to demonstrate the effectiveness of our theoretical results. 
\end{abstract}


\section{Introduction}

In multi-agent distributed optimization, the goal of agents is to cooperatively minimize the global objective function formed
	by the sum of local functions \cite{4749425}. Along with the penetration of multi-agent networks \cite{wang2007new,wang2009controllability,5404832}, distributed optimization
	has received ever-increasing attention in recent years \cite{4749425,6930814,7518617}. This is due to its wide practical applications in many
	areas such as large-scale sensor networks 
 \cite{7973152},
 and distributed power systems \cite{9141512}.

 In practical applications, optimization problems usually
 occur in uncertain environments. In practical applications, the accurate gradient is hard to obtain due to
uncertainties in communications and
environments. For instance, when using sensor networks to monitor a remote objective, e.g., bird migration patterns \cite{4497285}, uncertainty is inevitable due to the measurement noise \cite{nedic2018network}. Thus, it is necessary to study the
distributed optimization problem with stochastic gradients. In fact, various online gradient-free algorithms have been achieved.
	In \cite{Lee2017StochasticDA}, two decentralized stochastic variants of the Nesterov dual averaging method are proposed.
	In \cite{10019585}, an online
	gradient-free distributed projected gradient descent algorithm is presented via a two-point policy.
	In \cite{Xiong2022DistributedOO}, an online distributed constrained optimization problem under time-varying unbalanced digraphs without an explicit subgradient is studied.

  It is worth noting that all the aforementioned investigations
 study the the bounds of the regrets in expectation. For the algorithm with sublinear regret in expectation, it is necessary to run an algorithm in large numbers of rounds to eliminate the error between the unbiased estimator and the accurate gradient. However, in some practical problems, we can only run the algorithm for a few rounds or even once, so achieving the high probability bound of the regret is desired. For instance, in the problem of tracking a target \cite{Shahrampour2016DistributedOO}, agents need to track the target as fast as possible. In fact, in centralized optimization, the high probability convergence has already been studied \cite{madden2021highprobability,hong2023high,pmlr-v202-liu23aa}. Unfortunately, the  results in \cite{madden2021highprobability,hong2023high,pmlr-v202-liu23aa} are not applicable to the distributed cases. More recently, high probability convergence is studied in distributed settings. In \cite{10295561}, the high probability convergence of a distributed stochastic gradient algorithm are provided. In \cite{Yang2024}, the high probability bound of the dynamic regret is derived via a  distributed stochastic mirror descent.
 
 All works \cite{madden2021highprobability,hong2023high,pmlr-v202-liu23aa,10295561,Yang2024} above focus on the light-tailed noise. 
For light-tailed noise, such as the sub-Gaussian noise, large variety of concentration
 techniques are applicable since its tail decays faster than exponential distribution \cite{10295561}. Even though the light-tailed noise assumption is intuitive, in domains like evolutionary ecology \cite{Jourdain2012LvyFI}, audio source separation \cite{7177973}, finance \cite{1997Fractals} and machine learning \cite{pmlr-v97-simsekli19a}, the assumption seems invalid, which implies the necessity of studying heavy-tailed setting.

In this paper, an online distributed optimization problem is studied via a network of agents, where each agent only has access to a stochastic gradient of its own objective function, and can communicate with its neighbors via a network. To address this problem, an online distributed clipped stochastic gradient descent (ODCSGD) algorithm is proposed, and dynamic regrets are used to capture the performance of the algorithm. Furthermore, the high probability bounds of dynamic regrets are analyzed for both convex and non-convex objective functions. The main contributions are summarized as follows.

1. Here we consider the scenario where the stochastic gradient satisfy the heavy-tailed noise condition, as opposed to the lighted-tailed noise condition  in \cite{madden2021highprobability,hong2023high,pmlr-v202-liu23aa,10295561,Yang2024} . More specifically, if a noise follows a heavy-tailed distribution, then the variance of the noise is not necessarily bounded. Thus, the heavy-tailed noise condition is mathematically more general than the lighted-tailed noise condition. 

2. Different from the existing results on bounds of the regrets in expectation \cite{Lee2017StochasticDA, 10019585,Xiong2022DistributedOO,Shahrampour2016DistributedOO}, we study the high probability bound of the regret. To ensure the high probability bound of the stochastic gradient, the clippings strategy is employed. Combining the probability theory, convex optimization and consensus theory, we prove that if the graph sequence is $B$-strongly connected, and the objective function is convex, then a sublinear high probability bound of dynamic regret is guaranteed. For non-convex settings, if the variation of the objective function grows within a certain rate, a sublinear high probability bound of dynamic regret is obtained. 

To the best of our knowledge, this paper is the first one to study the online distributed optimization with stochastic gradient satisfying the heavy-tailed noise condition. Our approaches and results guarantee sublinear regret bounds with only
a single run of the algorithm. Thus, compared with results on the sublinear bounds of regrets in expectation \cite{Lee2017StochasticDA, 10019585,Xiong2022DistributedOO,Shahrampour2016DistributedOO}, our results are, mathematically and practically, more efficient and rigorous, and are applicable to wider areas such as distributed learning systems and distributed tracking problems where it is desired to run algorithms in a few rounds.

\textbf{Notations}: 
 $[n]$ represents the set $\{1,2,\cdots,n\}$ for any integer $n$. $\langle\cdot, \cdot\rangle$ is the standard inner product operator. $\mathbb{E}[\cdot]$ is the expectation operator. $\|\cdot\|$ represents the $2$-norm operator.

  This paper is organized as follows. In Section \ref{Sec2}, mathematical preliminaries on heavy-tailed noise and graph theory are introduced. In Section \ref{Sec3}, we
 state our main result and give its proof. In Section \ref{Sec4}, simulation
 examples are presented. Section \ref{Sec5} concludes the whole paper.
 \section{Preliminaries and Problem Formation}\label{Sec2}
 \subsection{Online Distributed Optimization with Noisy Gradient Information}
Let us describe a scenario of online distributed optimization. Consider a multi-agent system consisting of  $N$  agents, labeled by set  $\mathcal{V}=[N]$. Agents communicate with each other via a time-varying graph sequence  $\mathcal{G}_t$. For agent  $i \in \mathcal{V}$, a set of cost functions are given by  $\left\{f_{i,1}, \cdots, f_{i,T}\right\}$
 , where  $f_{i,t}: \mathbb{R}^d \rightarrow \mathbb{R}$ is not necessarily convex for any  $t \in[T],~ T \in \mathbb{N}$ is unknown to the agents. At each iteration time  $t \in[T]$, agent $i$ selects a state $ x_{i,t} \in \mathbb{R}^d$. After the state is selected, the information associated with the local cost function  $f_{i,t}$  is received by agent  $i$, that is, information on cost functions is not available before decisions are made by agents. 
 The goal of agents is to cooperatively solve the following optimization problem:
	\begin{equation}\label{onlineopt}
		\begin{split}
			\min f_{t}(x)=\sum_{i=1}^{N} f_{i,t}(x), \text { subject to } x \in \mathbb{R}^d.
		\end{split}
	\end{equation}
	To evaluate the performance of the online algorithm, a measurement called dynamic regret is adapted, which can be defined as
	\begin{equation}\label{dynamicregret}
		\begin{split}
			\boldsymbol{REG}_{T}^{d}=\sum_{i=1}^{N}\sum_{t=1}^{T} f_{t}\left(x_{i,t}\right)-\sum_{t=1}^{T} f_{t}\left(x^{*}_t\right).
		\end{split}
	\end{equation}
	where $x^{*}_t=\operatorname{argmin} f_t(x)$. 
 
	For general non-convex optimization problems, seeking the minimizer is usually NP-hard \cite{10025380}. Here the offline goal is to find the stationary point each time. Accordingly, the dynamic regret is defined as \cite{Li2023OnlineDS}
	\begin{equation}\label{dynamicregret_non-convex}
		\begin{split}
			\boldsymbol{NREG}_{T}^{d}=\sum_{t=1}^{T}\sum_{i=1}^{N}\|\nabla f_{t}(x_{i,t})\|^2.
		\end{split}
	\end{equation}
 
 Achieving a sublinear bound for dynamic regret (\ref{dynamicregret}) and (\ref{dynamicregret_non-convex}) is rather challenging in the worst case when the objective functions change fast. However, the path length can be used to quantify the difficulty \cite{9117183}. The path length is defined as follows
	\begin{equation}
	C_T=\sum^{T}_{t=2} \|x_t^*-x_{t-1}^*  \| .
	\end{equation}
We employ the variation of the function
sequence to reflect the difficulty in the worst case \cite{doi:10.1287/opre.2015.1408}: 
\begin{align}\label{non-convex_path}
    D_T=\sum_{t=2}^{T}\|f_{t}-f_{t-1}\|_{\sup}
\end{align}
where $\|f_{t}-f_{t-1}\|_{\sup}:=\sup_{x\in \mathbb{R}^{m}}\|f_{t}(x)-f_{t-1}(x)\|$.
 
	Here we assume that agents only have access to a noisy gradient of $f_{i,t}$, denoted by $\widehat{\nabla}f_{i, t}: \mathbb{R}^d\to  \mathbb{R}^d$. For the estimation of the gradient, we make the following assumptions:
\begin{assu}\label{assu2}
		For any $x\in \mathbb{R}^d$,
		\begin{equation*}
			\mathbb{E}[\widehat{\nabla}f_{i,t}(x)|x]=\nabla f_{i,t}(x),~\forall~i\in [n].
		\end{equation*}
	\end{assu}
	\begin{assu}\label{assu3}
		For any $x\in \mathbb{R}^d$, $i\in [N]$, 
		\begin{align*}
		   \mathbb{E}[\|\widehat{\nabla}f_{i,t}(x)-\nabla f_{i,t}(x) \|^p|x]\le \sigma^p 
		\end{align*}
		for some $\sigma>0$ and $p\in (1,2]$.
	\end{assu}
 
 This is commonly referred to as heavy-tailed noise, as opposed to light-tailed noise such as those that are distributed according to sub-Gaussian distributions. Mathematically, heavy-tailed noise is more general in stochastic optimization. The  variance of the light-tailed noise, e.g. sub-Gaussian noise, is finite since its tail decays faster than Gaussian distribution. Thus, light-tailed noise satisfies Assumption \ref{assu3}, but not vice versa.
	\begin{assu}\label{assu4}
		For any $i\in [N]$, $f_{i,t}$ is $L$-smooth, i.e., its gradient is $L$-Lipchitz. For all $x,y\in \mathbb{R}^d$,
		\begin{equation*}
			|\nabla f_{i,t}(x)- \nabla f_{i,t}(y)|\leq L \| x-y\|.
		\end{equation*}
   Moreover, if $f_{i,t}$ is $L$-smooth, it has a quadratic upperbound, i.e., For all $x,y\in \mathbb{R}^d$,
        \begin{equation*}
		 f_{i,t}(y)\le f_{i,t}(x) +\langle \nabla f_{i,t}(x), y-x \rangle +\frac{L}{2}\|y-x\|^2.
		\end{equation*}
	\end{assu}

 \begin{assu}\label{assu5}
    For any $i\in [N]$ and $t\in [T]$, $x_{i,t}$ is almost surely bounded. That is, there exists $B_X>0$, such that $\|x_{i,t}\|\le B_X$. Moreover, there exists $B_g> 0$ such that $\|\nabla f_{i,t}(x_{i,t})\|\le B_g$ and $\|\nabla f_{i,t}(\bar{x}_t)\|\le B_g$. $\bar{x}_t$ is the average of $x_{i,t}$ for $i\in [N]$.
 \end{assu}
 
Agents communicate with each other via an directed graph $\mathcal{G}_t=(\mathcal{V}, W_t, \mathcal{E}_t)$, where $\mathcal{V}$ denotes the set of agents, $\mathcal{E}_t$ denotes the set of edges where the elements are denoted as $(i,j)$ if and only if agent $i$ can
receive a message from agent $j$ at time $t$, and $W_t=([W_t]_{ij})_{n\times n}$ denotes the weighted matrix. We denote the set of
incoming neighbors of agent $i$ at time $t$ by
\begin{equation*}
    \mathcal{N}_{i,t}=\{j|(i,j)\in \mathcal{E}_t \}\cup \{i\}.
\end{equation*}
For a fixed topology $\mathcal{G}_t=(\mathcal{V}, W, \mathcal{E})$, a path of length $r$ from $i_1$ to $i_{r+1}$ is a sequence of $r+1$ distinct nodes $i_1, \cdots, i_{r+1}$ such that $(i_{q},i_{q+1})\in \mathcal{E}$, for $q\in [r]$. If there exists a path between any two nodes, then $\{\mathcal{G}_t\}$ is said to be strongly connected. For $\{\mathcal{G}_t\}$, an $B$-edge set is defined as $\mathcal{E}_{B,t}=\bigcup_{l=0,\cdots,B-1}\mathcal{E}_{t+l}$ for some constant $B>0$. We call that $\{\mathcal{G}_t\}$ is $B$-strongly connected if the directed graph with vertex $\mathcal{V}$ and edge $\mathcal{E}_{B,t}$ is strongly connected for any $t\ge 0$.
	\begin{assu}\label{assu6}
		For all $t \ge 1$, the weighted graphs $\mathcal{G}_t=(\mathcal{V}, W_t, \mathcal{E}_t)$ satisfy:\\
  (a) There exists a scalar $\eta \in (0, 1)$ such that $[W_t]_{ij}\ge \eta  $ if $j\in \mathcal{N}_{i,t}$. Otherwise, $[W_t]_{ij}=0 $.\\
  (b) The weighted matrix is doubly stochastic, i.e., $\sum_{i=1}^{n}[W_t]_{ij}=\sum_{j=1}^{n}[W_t]_{ij}=1$.\\
  (c) $\{\mathcal{G}_t\}$ is $B$-strongly connected. That is, there exists a scalar $B > 0$ such that the graph
$(\mathcal{V},\mathcal{E}_{B,t})$ is strongly connected for any $t \ge 1$.
	\end{assu}
\subsection{Online Distributed Clipped Stochastic Gradient Descent}
 To solve problem (\ref{onlineopt}), we propose the following ODCSGD algorithm
\begin{equation}\label{algorithmonline}
\begin{cases}
    y_{i,t}=\sum_{j=1}^{N}\left[W_{t}\right]_{i j}x_{j,t}\\
    \widetilde{\nabla} f_{i,t}(x_{i,t})=\text{clip}_{\lambda_t}(\widehat{\nabla} f_{i,t}(x_{i,t}))\\
    x_{i,t+1}=y_{i,t}-\eta_t\widetilde{\nabla} f_{i,t}(x_{i,t}).
\end{cases}
\end{equation}
where $\{\eta_{t}\}_{t=1}^{T}$ is the global step size sequence, $x_{i,t}\in \mathbb{R}^d$ represents the state of agent $i$ with initial state $x_{i,1}\in \mathbb{R}^d$. $\widehat{\nabla}f_{i,t}: \mathbb{R}^d\to \mathbb{R}^d$ is the non-biased estimation of the gradient ${\nabla}f_{i,t}$. The step size sequence is non increasing and positive and the clipping operator $\text{clip}_{\lambda_t}(\cdot): \mathbb{R}^d\to \mathbb{R}^d$ is 
 \begin{align}\label{clipping_operator}
     \text{clip}_{\lambda_t}(y)=\min \{1,\frac{\lambda_t}{\|y\|}\}y.
 \end{align}
\section{Main Results}\label{Sec3}
Before presenting our main results, we provide some useful lemmas. First we analyze the network error, i.e., the error between each agent’s state
 and their average value at each iteration under (\ref{algorithmonline}). 
\begin{lemma}[\cite{sundhar2012new}]\label{lem1}
    Let Assumption \ref{assu6} holds for a sequence of weight matrices  $\left\{W_{t}\right\}_{t \geq 1}$. Consider a set of sequences  $\left\{\boldsymbol{\xi}_{i, t}\right\} $ for $ i \in [N]$ defined by the following relation:
    \begin{align*}
        \boldsymbol{\xi}_{i, t+1}=\sum_{j=1}^{N}\left[W_{t}\right]_{i j} \boldsymbol{\xi}_{j, t}+\boldsymbol{\epsilon}_{i, t+1}, \text { for } t \geq 1
    \end{align*}
Let $ \overline{\boldsymbol{\xi}}_{t} $ denote the average of  $\boldsymbol{\xi}_{i, t} $ for $ i \in[N] $, i.e.,  $\overline{\boldsymbol{\xi}}_{t}=   \frac{1}{N} \sum_{i=1}^{N} \theta_{i, t} $. Then,
\begin{align*}
    &\left\|\boldsymbol{\xi}_{i, t+1}-\overline{\boldsymbol{\xi}}_{t+1}\right\| \leq N \gamma \beta^{t} \max _{j}\left\|\boldsymbol{\xi}_{j, 1}\right\| \\
	&+\gamma \sum_{\ell=1}^{t-1} \beta^{t-\ell} \sum_{j=1}^{N}\left\|\boldsymbol{\epsilon}_{j, \ell+1}\right\|+\frac{1}{N} \sum_{j=1}^{N}\left\|\boldsymbol{\epsilon}_{j, t+1}\right\|+\left\|\boldsymbol{\epsilon}_{i, t+1}\right\|
\end{align*}
where $ \gamma $ and  $\beta $ are defined as
\begin{align}\label{gamma_eta}
    \gamma=\left(1-\frac{\eta}{2 N^{2}}\right)^{-2} \quad \beta=\left(1-\frac{\eta}{2 N^{2}}\right)^{\frac{1}{B}}.
\end{align}
\end{lemma}$\\$
Note that $x_{i,t}$ can be represented by the iteration in Lemma \ref{lem1}. By (\ref{gamma_eta}) and algorithm (\ref{algorithmonline}), we can bound the network error as follows
\begin{lemma}\label{lem_networkerr}
    Under Assumption \ref{assu6}, by algorithm (\ref{algorithmonline}), for all $ i \in [N] $ and $ t \geq 1$, 
\begin{align*}
	\left\|x_{i, t+1}-\bar{x}_{t+1}\right\|& \leq  N\gamma \beta^{t} R_1+2 \lambda_{t}\eta_t+
 N\gamma \sum_{l=1}^{t-1}\beta^{t-l}\lambda_{l}\eta_l,
\end{align*}
where $\bar{x}_t=\frac{1}{N}\sum_{i=1}^{N}x_{i,t} $, $ \gamma $ and $ \beta $ are defined in (\ref{gamma_eta}) and $R_1=\max_{i}\|x_{i,1}\|$.
\end{lemma}

\textit{Proof:} By algorithm (\ref{algorithmonline}), it follows that
\begin{align*}
    x_{i,t+1}=\sum_{j=1}^{N}[W_t]_{ij}x_{j,t}-\eta_t\widetilde{\nabla}f_{i}(x_{i,t}).
\end{align*}
By the second iteration in algorithm (\ref{algorithmonline}) and (\ref{clipping_operator}), we have $\|\eta_t\widetilde{\nabla}f_{i}(x_{i,t})\|\le \eta_t \lambda_t$, by Lemma \ref{lem_networkerr}, and the decreasing property of $\eta_t\lambda_t$, the lemma is proved.

Moreover, we have
\begin{lemma}\label{cor_network^2}
    Under Assumption \ref{assu6}, by algorithm (\ref{algorithmonline}), for all $ i \in [N] $ and $ t \geq 1$, 
\begin{align*}
	&\sum_{t=1}^{T}\left\|x_{i, t}-\bar{x}_{t}\right\|^2 \leq  3N^2\gamma^2R_1^2\frac{1}{1-\beta^2}+12\sum_{t=1}^{T}\eta_t^2\lambda_t^2\\
 &+3N^2\gamma^2\frac{1}{(1-\beta)^2}\sum_{t=1}^{T}\eta_t^3\lambda_t^2,
\end{align*}
\end{lemma}

\textit{Proof:} See Appendix \ref{proof_network^2}.

\subsection{Convex Case}
When the objective functions are convex, the following lemma presents the bound for the term $\eta_t (f_{i, t}\left(x_{i, t}\right)-f_{i, t}\left(x_t^{\star}\right ))$. 
\begin{lemma}\label{thm_1}
    Under Assumptions \ref{assu5}-\ref{assu6}, by algorithm (\ref{algorithmonline}), 
\begin{equation}\label{thm_main_error}
\begin{aligned}
    &\eta_t (f_{i, t}\left(x_{i, t}\right)-f_{i, t}\left(x_t^{\star}\right ))-\eta_t\langle \theta_{i,t},x_{i, t}-x_t^{\star}    \rangle \\&\le \frac{1}{2} (\|x_{i, t}-x_{t}^{\star}\|^2-\|x_{i, t+1}-x_t^{\star}\|^2)+5\eta_t^2\lambda_t^2 \\
    &+2\lambda_t\eta_t(N\gamma \beta^{t-1} R_1+
 N\gamma \sum_{l=1}^{t-1}\beta^{t-l}\lambda_{l}\eta_l)
\end{aligned}
\end{equation}
where $\langle \theta_{i,t},x_{i, t}-x^{\star}\rangle=\langle\widetilde\nabla f_{i,t}(x_{i,t})-\nabla f_{i,t}(x_{i,t}),x_{i, t}-x_t^{\star}    \rangle$.
\end{lemma} 

\textit{Proof:} See Appendix \ref{Proofthm1}. 

Due to the fact that $\theta_{i,t}$ is not unbiased, it is difficult to analyze the high probability bound of the term $\eta_t\langle \theta_{i,t},x_{i, t}-x_t^{\star}    \rangle$. Let
\begin{align*}
    &\theta_{i,t}^u=\widetilde{\nabla}f_{i,t}(x_{i,t})-\mathbb{E}[\widetilde{\nabla}f_{i,t}(x_{i,t})|\mathcal{F}_t],\\&
    \theta_{i,t}^b=\mathbb{E}[\widetilde{\nabla}f_{i,t}(x_{i,t})|\mathcal{F}_t]-\nabla f_{i,t}(x_{i,t}).
\end{align*}
where $\mathcal{F}_t=\sigma(\widehat{\nabla} f_{1},\cdots,\widehat{\nabla} f_{t-1})$ denotes the $\sigma$-field generated by the unbiased estimator by $t$ ($(\mathcal{F}_t)_{t\ge 0}$ is also known as the natural filtration).
Based on results in \cite{nguyen2023high}, preliminary bounds of $\theta_{i,t}^u$  and $\theta_{i,t}^b$ are presented.
\begin{lemma}[\cite{nguyen2023high}]\label{lem_heavytaied}
    For $t\ge1$ and $i\in [N]$, under Assumption \ref{assu2}, we have
    \begin{align*}
        \|\theta_{i,t}^u\|\le 2\lambda_t.
    \end{align*}
    Furthermore, if $\|\nabla f_{i,t}(x_{i,t})\|\le \frac{\lambda_t}{2}$, then
    \begin{align*}
        \|\theta_{i,t}^b\|\le 4\sigma^p\lambda_{t}^{1-p},~~~\mathbb{E}[\|\theta_{i,t}^u\|^2|\mathcal{F}_t]\le 16\sigma^p\lambda_t^{2-p}.
    \end{align*}
\end{lemma}

Then, the high probability bound of the term $\eta_t\langle \theta_{i,t},x_{i, t}-x_t^{\star} \rangle$ is provided as follows.
\begin{lemma}\label{lem_highproboundonline}
    Under Assumptions \ref{assu2}-\ref{assu6}, for $\lambda_t=2B_g t^{\alpha}$ and $\eta_t=\frac{1}{(at+b)^{\kappa}}$, $\kappa>2\alpha>0$, by algorithm (\ref{algorithmonline}), for any $\delta\in (0,1)$, with probability at least $1-\delta$,
    \begin{align*}
    &\sum_{t=1}^{T} \eta_t\langle \theta_{i,t},x_{i, t}-x_t^{\star} \rangle\\
    &\le 8 B_X \sum_{t=1}^{T} \eta_t \lambda_t (\frac{\sigma}{\lambda_t})^p+ \frac{16}{3} B_XB_g \log \frac{2}{\delta}+ \sqrt{2F}.
    \end{align*}
\end{lemma}

\textit{Proof:} See Appendix \ref{proofhighbound}.

Now we are in a position to present our main result.
\begin{thm}\label{thm_dynamicregret}
Under Assumptions \ref{assu2}-\ref{assu6}, for $\lambda_t=2B_g t^{\alpha}$ and $\eta_t=\frac{1}{(at+b)^{\kappa}}$, $\kappa>2\alpha>0$, by algorithm (\ref{algorithmonline}), for any $\delta\in (0,1)$, with probability at least $1-\delta$,
 \begin{align*}
    \boldsymbol{REG}_{T}^{d}\le& P\frac{1}{\eta_T}+Q\sum_{t=1}^{T}\frac{\eta_t^2\lambda_t^2}{\eta_T}+R\sum_{t=1}^{T}\frac{\eta_t^2\lambda_t}{\eta_T}\\
    &+2B_XC_T\frac{1}{\eta_T}+\frac{16N}{3} B_XB_g \log \frac{2}{\delta}\frac{1}{\eta_T}\\
    &+8N B_X \sigma^p \sum_{t=1}^{T}  \lambda_t^{1-p} \
    \frac{\eta_t}{\eta_T} \\
    &+8N\frac{\sqrt{\sigma^p }B_X}{\eta_T}\sqrt{\sum_{t=1}^{T}\lambda_t^{2-p}\eta_t^2}
\end{align*}
where
$P=(2B_X \lambda_1+2B_gN^2)N\gamma\eta_1\frac{1}{1-\beta}+\frac{1}{2}\sum_{i=1}^{N}\|x_{i,1}-x_{1}^*\|^2$, $Q=(5+2N\gamma\frac{1}{1-\beta})N$
and $R=2B_gN^2(2+N\gamma\frac{1}{1-\beta})$.
\end{thm}

\textit{Proof:} See Appendix \ref{proofstatic}.

Based on Theorem \ref{thm_dynamicregret}, we have the following corollary.
\begin{corollary}
    Under the same conditions
stated in Theorem \ref{thm_dynamicregret}, if $\lambda_t=2B_g t^{\alpha}$ and $\eta_t=\frac{1}{(at+b)^{\kappa}}$ with $\kappa=\frac{p+1}{2p}$ and $\alpha=\frac{1}{2p}$, then with probability at least $1-\delta$,
\begin{align}\label{convex_dynamicregret}
   \boldsymbol{REG}_{T}^{d}\sim \mathcal{O}(T^{\frac{1+p}{2p}}(1+C_T+\log \frac{1}{\delta})).
  \end{align}
\end{corollary}

\textit{Proof:} Note that for $a<1$, 
\begin{align*}
    \sum_{t=1}^{T}\frac{1}{t^a}\le \int_{0}^{T} \frac{1}{s^a} d t= \frac{T^{1-a}}{1-a}
\end{align*}
which yields that
\begin{align*}
    &\sum_{t=1}^{T}\frac{\eta_t^2\lambda_t^2}{\eta_T}\le \frac{T^{2\alpha-\kappa+1}}{2\alpha-\kappa+1},
    \sum_{t=1}^{T}\frac{\eta_t^2\lambda_t}{\eta_T}\le \frac{T^{\alpha-\kappa+1}}{\alpha-\kappa+1},\\
    &\sum_{t=1}^{T}  \lambda_t^{1-p} \frac{\eta_t}{\eta_T}\le \frac{T^{\alpha(1-p)+1}}{\alpha(1-p)+1},\\
    &\sum_{t=1}^{T}\lambda_t^{2-p}\eta_t^2\le
    \frac{T^{(2-p)\alpha-2\kappa+1}}{(2-p)\alpha-2\kappa+1}.
\end{align*}
Note that when $\kappa=\frac{p+1}{2p}$ and $\alpha=\frac{1}{2p}$, $\max \{\kappa, 2\alpha-\kappa+1,\alpha(1-p)+1, \frac{1}{2}+\alpha(1-\frac{p}{2})\}=\frac{p+1}{2p}$. The validity of the corollary is proven.
\begin{rem}
    From Corollary 1, the sublinearity of the bound in (\ref{convex_dynamicregret}) is influenced by term  $T^{\frac{1+p}{2p}} \log \frac{1}{\delta}$ . Note that the value of  $\log \frac{1}{\delta} $ increases slowly as the value of failure probability  $\delta$  decreases. For example, due to the facts that  $\ln 10^{2}=4.61$, $\ln 10^{3}=6.91$, the term  $T^{\frac{1+p}{2p}} \log \frac{1}{\delta} $ sublinearly increases as  $4.61 T^{\frac{1+p}{2p}}$ ,  $6.91 T^{\frac{1+p}{2p}}$, with probabilities at least  $99.99 \%$, $99.999 \%$, respectively. Hence, the sublinearity of term  $T^{\frac{1+p}{2p}} \log \frac{1}{\delta}$  with a probability close to  $100 \%$  can be ensured by running algorithm (\ref{algorithmonline}) in a single round. Moreover, the sublinearity of the bound in (\ref{convex_dynamicregret}) is also influenced by  $C_T$. If  $C_T$  is sublinear with  $T^{\frac{-1+p}{2p}}$ , i.e.,  $\lim _{T \rightarrow \infty} \frac{C_T}{T^{\frac{-1+p}{2p}}}=0$, then  $\boldsymbol{REG}_{T}^{d}$ has a sublinear bound with high probability. If the minimizer changes fast, $C_T$ may be linear with $T^{\frac{-1+p}{2p}}$, then the sublinearity of the regret can not be guaranteed.
\end{rem}

\subsection{Non-convex Case}
The result on the non-convex case is provided as follows.

\begin{thm}\label{thm_non-convex}
Under Assumptions \ref{assu2}-\ref{assu6}, for $\lambda_t=2B_g t^{\alpha}$ and $\eta_t=\frac{1}{(at+b)^{\kappa}}$, $\kappa>2\alpha>0$, by algorithm (\ref{algorithmonline}), for any $\delta\in (0,1)$, with probability at least $1-\delta$,
 \begin{align*}
    \frac{\boldsymbol{NREG}_{T}^{d}}{2N}\le& P \frac{1}{\eta_{T}}+Q\sum_{t=1}^{T}\frac{n_{t}^{2} \lambda_{t}}{\eta_{T}}+R \sum_{t=1}^{T} \eta_{t}^{2} \lambda_{t}^{2}\\
    &+4 N B_{g} \sigma^{p} \sum_{t=1}^{T}\lambda_{t}^{1-p} \frac{\eta_{t}}{\eta_{T}}+N D_T \frac{1}{\eta_{T}}\\
    &+4 N B_{g} \frac{\sigma^{\frac{p}{2}}}{\eta_{T}} \sqrt{\sum_{t=1}^{T} \lambda_{t}^{2-p} \eta_{t}^{2}}+\frac{8}{3} B_{g}^{2} \log \frac{2}{\delta} \frac{1}{\eta_{T}} \\
    &+ 3N^3L^2\gamma^2\frac{1}{(1-\beta)^2}\sum_{t=1}^{T} \eta_{t}^{3} \lambda_{t}^{2}
\end{align*}
where
$P=N\left(f_{1}\left(\bar{x}_{1}\right)-f_{T+1}\left(\bar{x}_{T+1}\right)\right)+N \gamma \frac{\eta_{1}}{1-\beta} R_{1} B_{g} L+36N^3L^2R_1^2\frac{1}{1-\beta^2} $, $Q=BgL\left(2+N \frac{\gamma}{1-\beta}\right)$ and $R=(\frac{L}{2}+12NL^2)$.
\end{thm}

\textit{Proof:} See Appendix \ref{proof_non-convex}.

Moreover, the high probability bound of the regret via some specific step sizes can be obtained.
\begin{corollary}
    Under the same conditions
stated in Theorem \ref{thm_non-convex}, if $\lambda_t=2B_g t^{\alpha}$ and $\eta_t=\frac{1}{(at+b)^{\kappa}}$ with $\kappa=\frac{p+1}{2p}$ and $\alpha=\frac{1}{2p}$, then with probability at least $1-\delta$,
\begin{align*}
   \boldsymbol{NREG}_{T}^{d}\sim \mathcal{O}(T^{\frac{1+p}{2p}}(1+D_T+\log \frac{1}{\delta})).
  \end{align*}
\end{corollary}
\begin{rem}
Different from convex settings, $D_T$ is employed as a complexity measure of the problem environment. If  $D_T$  is sublinear with  $T^{\frac{-1+p}{2p}}$, then  $\boldsymbol{NREG}_{T}^{d}$  has a sublinear bound with high probability. Accordingly, algorithm (\ref{algorithmonline}) performs well.
\end{rem}
\section{Numerical Simulations}\label{Sec4}
\begin{figure}
     \centering
     \includegraphics[scale=0.15]{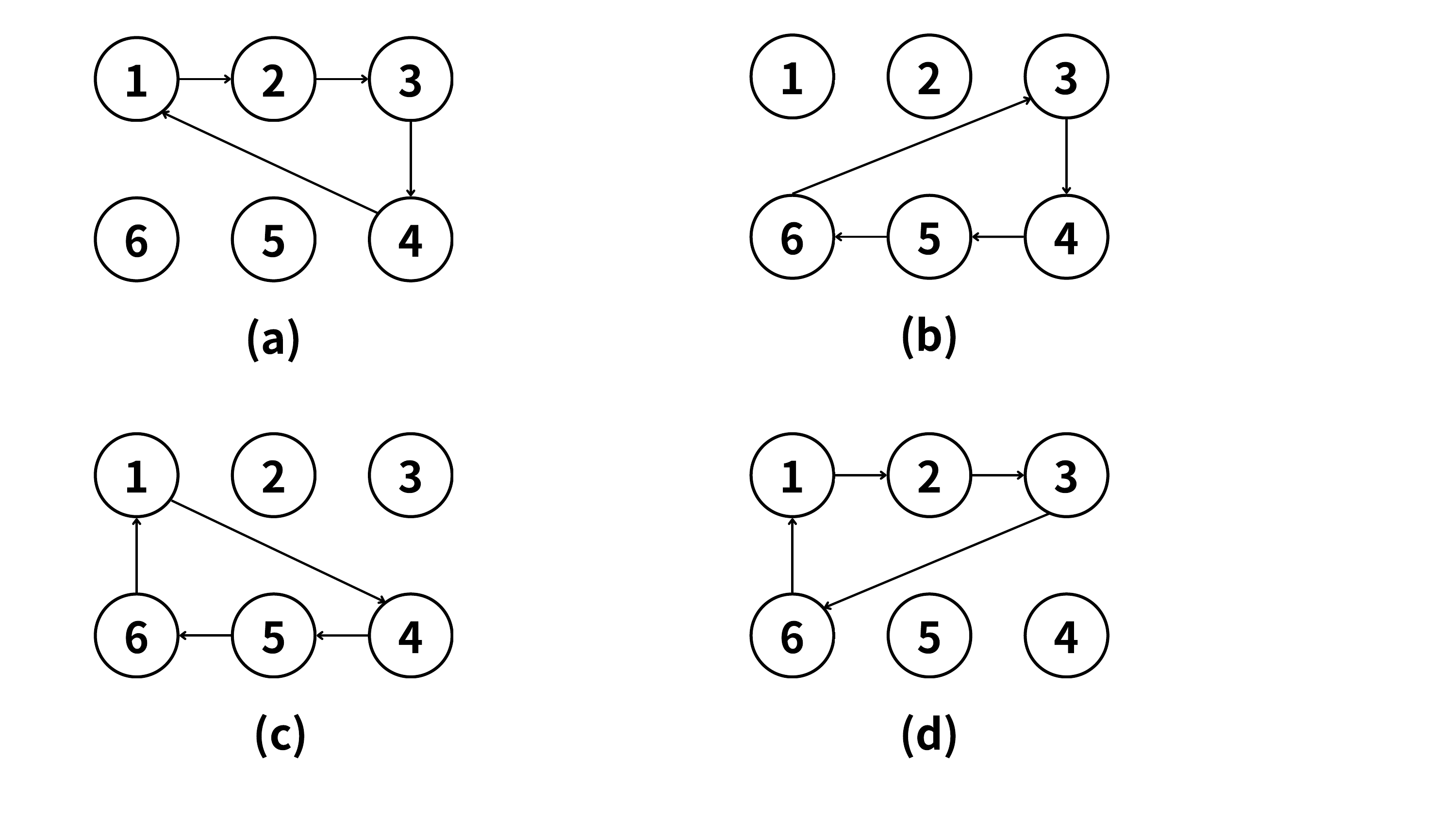}
     \caption{The Time-varying Graph}
     \label{fig:time-varying}
 \end{figure}

 In this section, we use our methods to solve a distributed
 tracking problem, where sensors aim to cooperatively track
 a target. Consider a sensor network consisting of 6 sensors,
 labeled by set $V = \{1,\cdots ,6\}$. 
	Each sensor communicates
 with its neighbors via a time-varying graph shown in Fig.\ref{fig:time-varying}.
 
 The weight
 of each edge in Fig.\ref{fig:time-varying} is assumed to be $0.8$. The switching order is given by $(a)\to (b)\to (c)\to (d)\to (a)\to \cdots$. Note that
 the union of the possible graphs is strongly connected. Then,
 the connectivity of the graph in Fig.\ref{fig:time-varying} satisfies conditions in
 Assumption \ref{assu6} with $B = 4$.

 The position of the target, denoted by $x_t^*\in \mathbb{R}^2$, evolves as
	\begin{equation*}
		\begin{cases}
			x_{t+1}^*=x_t^*+\omega_t+v_t, \\
			x_0^*=(10,10)^\top,\quad t=0,1,\cdots,T-1
		\end{cases}
	\end{equation*}
	where $\omega_t$ is the Gaussian noise with zero expectation and $(\frac{1}{T})^2$ variance. $v_t$ is the nonlinear path of the $x_t^*$:
	\begin{equation*}
	v_t=-10\begin{pmatrix}
			\frac{1}{T}, &
			\sqrt{\frac{t}{T} }-\sqrt{\frac{t-1}{T} }
		\end{pmatrix}^{\top}.	
	\end{equation*}
	Agents aim to track the position of the moving object collaboratively. At time $t$, each agent observes one of the coordinates of $x_t^*$ with $\xi_{i,t}$, i.e.,
	\begin{equation*}
		z_{i,t}=e_{k_i}^{\top}x_t^*
	\end{equation*}
	where $e_{k_i}$ is the $k$-th unit vector in the standard basis of $\mathbb{R}^2$ for $k \in \{1,2\}$. We divide the agents into two groups, where $k_i$ is the remainder of $i$ divided by $2$, $i=1,2,\cdots,6$.
 
 \subsection{Convex Case}
 We choose the average of the local square loss as the measure of the tracking error:
	\begin{equation}\label{simuprob}
		f(x)=\frac{1}{n} \sum^{n}_{i=1}f_{i}(x)=\frac{1}{n} \sum^{n}_{i=1} \mathbb{E}
		\left[(z_{i,t}-e_{k_i}^{\top}x_t)^2\right].
	\end{equation}
	 Note that the local gradient is not accurate, i.e.,
	\begin{equation*}
		\widehat{\nabla}f_{i, t}= (e_{k_i}^{\top}(x-x^*)+\xi_{i,t})e_{k_i}
	\end{equation*}
	where $\xi_{i,t}$ is the heavy-tailed noise. We consider a typical noise, the t-distribution $t_2$ with the probability density function
	\begin{align*}
	    t_2(x)=\frac{\Gamma(\frac{3}{2})}{\Gamma(\frac{1}{2})
	   \sqrt{2\pi}}(1+\frac{x^2}{2})^{-\frac{3}{2}}
	\end{align*}
	where $\Gamma(x)=\int_{0}^{
	+\infty}t^{x-1}e^{-t}dt$. We can verify that $t_2$ has zero expectation and unbounded variance. 
	Then we solve problem (\ref{simuprob}) by our algorithm. In simulation, we choose $T=5000$ with the parameter $\eta_t=\frac{1}{(0.5t+10)^{0.5}}$ and $\lambda_t=2t^{0.1}$.
 
\begin{figure}
        \centering
        \includegraphics[scale=0.4]{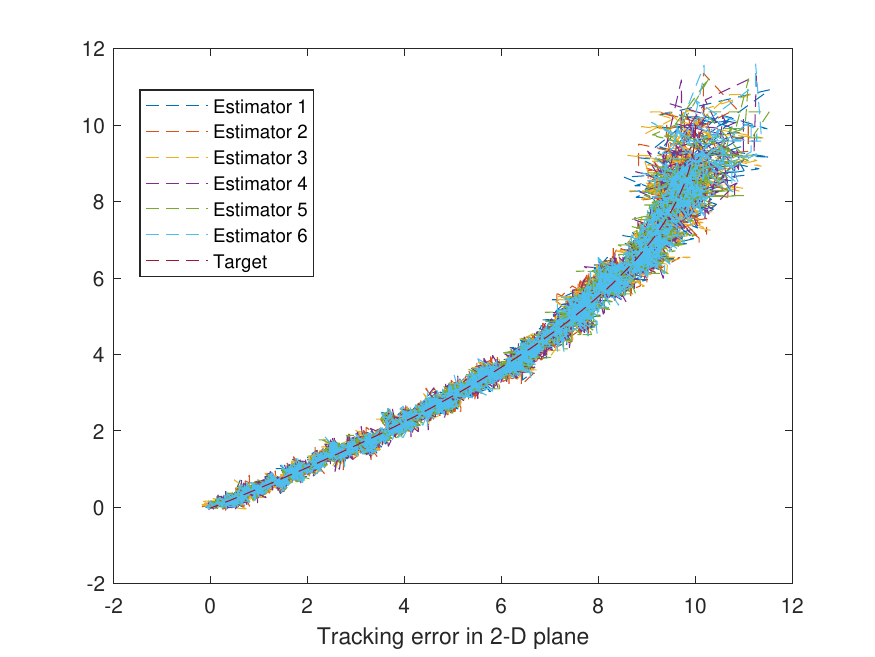}
        \caption{Tracking Error in 2-D Plane}
        \label{fig:Tracking error 2d}
    \end{figure}

\begin{figure}
        \centering
        \includegraphics[scale=0.4]{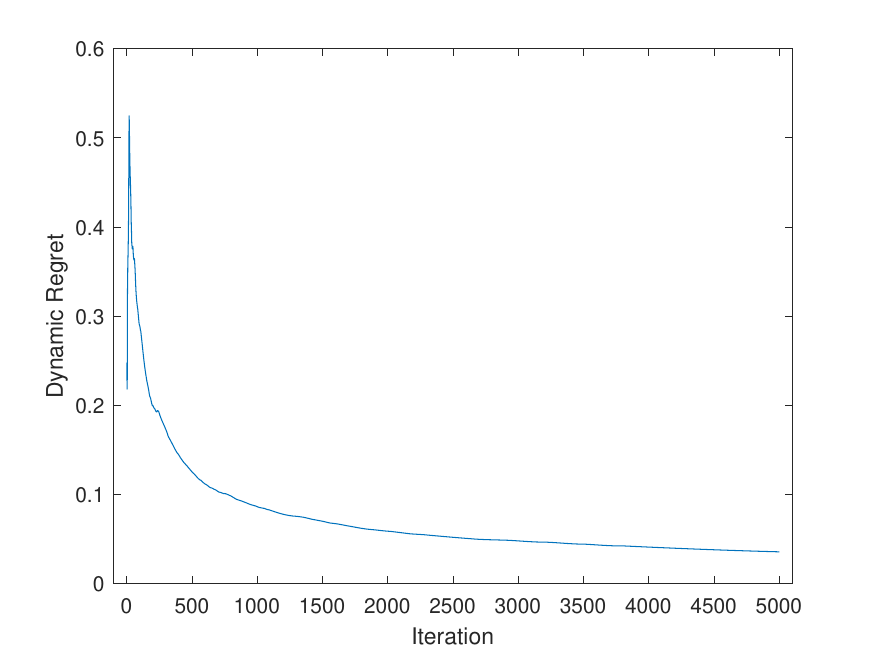}
        \caption{The Dynamic Regret Averaged by Iteration}
        \label{fig:dynamic_regret}
    \end{figure}

    The initial values are randomly chosen from $[9,10] \times [9,10]$. By running our algorithm in a single round, we present the estimation of the target's trajectory of each agent in Fig. \ref{fig:Tracking error 2d}. 
    Additionally, we plot the dynamic regret divided by time, as depicted in Fig. \ref{fig:dynamic_regret}. 
  \subsection{Non-convex Case}
\begin{figure}
     \centering
     \includegraphics[scale=0.4]{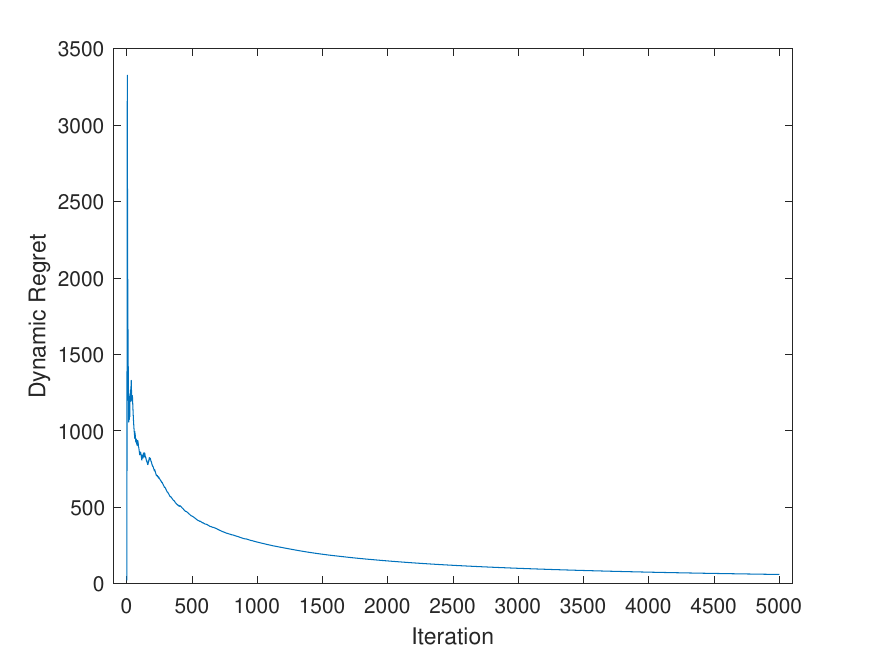}
     \caption{$\boldsymbol{NREG}_{T}^d/T$}
     \label{fig:nonconvex_dynamic_regret}
 \end{figure}

 In non-convex settings, the value function is modified as follows
 \begin{equation}\label{simuprob2}
		f(x)=\frac{1}{n} \sum^{n}_{i=1}f_{i}(x)=\frac{1}{4n} \sum^{n}_{i=1} \mathbb{E}
		\left[(z_{i,t}^2-e_{k_i}^{\top}x_t^2)^2\right].
	\end{equation}
 By algorithm (\ref{algorithmonline}), with $\eta_t=\frac{1}{(0.5t+10)^{0.4}}$ and $\lambda_t=2t^{0.1}$, we plot the corresponding $\boldsymbol{NREG}_{T}^{d}/T $ which can be referred to Fig. \ref{fig:nonconvex_dynamic_regret}.
 Note that both
	$\boldsymbol{REG}_{T}^{d}/T $ and $\boldsymbol{NREG}_{T}^d/T$ decay, indicating  $\boldsymbol{REG}_{T}^{d} $ and $\boldsymbol{NREG}_{T}^d$ grow sublinearly. These observations are consistent with our obtained results.
	
\section{Conclusion}\label{Sec5}
In this paper, we have proposed an ODCSGD algorithm for
 the distributed optimization problem.  
 By implementing the algorithm,
 every agent adjusts its state value by 
 the clipped gradient estimation from its local cost function and the local
 state information received from its immediate neighbors.
 We show that
if the time-varying graph sequence is $B$-strongly connected and the objective function is convex, the corresponding high probability bound of the dynamic regret is $\mathcal{O}(T^{\frac{1+p}{2p}}(1+C_T+\log \frac{1}{\delta}))$. 
Moreover, in the non-convex settings, the corresponding high probability bound of the regret is $\mathcal{O}(T^{\frac{1+p}{2p}}(1+D_T+\log \frac{1}{\delta}))$. Numerical
 simulations have been presented to demonstrate the
 effectiveness of our theoretical results. 
 
 When constraints exist in distributed optimization problems,
 more complicated distributed stochastic algorithms, such as the one based on the primal-dual strategy, are needed. 
 How to study the bound of regret of those distributed algorithms in high
 probability is still an interesting open topic. 
\section{Appendix}
\subsection{Proof of Lemma \ref{cor_network^2}}\label{proof_network^2}
\textit{Proof of Lemma \ref{cor_network^2}:} By Lemma \ref{lem_networkerr}, we obtain that
\begin{align*}
    &\|x_{i,t}-\bar{x}_t\|^2\le  3\Big(\Big.(N\gamma \beta^{t} R_1)^2+4 \lambda_{t}^2\eta_t^2\\
    &+
 N^2\gamma^2 (\sum_{l=1}^{t-1}\beta^{t-l}\lambda_{l}\eta_l)^2\Big)\Big. .
\end{align*}
The inequality results from $(a_1+\cdots+a_n)^2\le n(a^2+\cdots+a_n^2)$. By the Cauchy-Schwartz inequality, we can derive that
\begin{align*}
    &\sum_{t=1}^{T}\eta_t(\sum_{l=1}^{t-1}\beta^{t-l}\lambda_{l}\eta_l)^2\le \sum_{t=1}^{T}\eta_t(\sum_{l=1}^{t-1}\beta^{t-l})\sum_{l=1}^{t-1}\beta^{t-l}(\lambda_{l}\eta_l)^2\\
    &\le \frac{1}{1-\beta}\sum_{l=1}^{T-1}\sum_{t=l+1}^{T-1}\beta^{t-l}(\lambda_{l}\eta_l)^2\eta_t\le (\frac{1}{1-\beta})^2\sum_{l=1}^{T-1}\lambda_{l}^2\eta_l^3.
\end{align*}
Summing over $t\in[T]$, we derive the result.
\subsection{Proof of Lemma \ref{thm_1}}\label{Proofthm1}
\textit{Proof of Lemma \ref{thm_1}:} By the convexity of $f_{i,t}$, we derive that
\begin{align*}
	&f_{i, t}\left(x_{i, t}\right)-f_{i, t}\left(x_t^{\star}\right )\le\left \langle\nabla f_{i,t} (x_{i,t}),x_{i, t}-x_t^{\star}   \right \rangle\\
	&= \langle\widetilde\nabla f_{i,t} (x_{i,t}) ,x_{i, t}-x_t^{\star} \rangle\\
	&+\langle\widetilde\nabla f_{i,t} (x_{i,t})-\nabla f_{i,t} (x_{i,t}) ,x_{i, t}-x_t^{\star} \rangle\\
	&= \langle\widetilde\nabla f_{i,t} (x_{i,t}) ,x_{i, t+1}-x_t^{\star} \rangle+\langle\widetilde\nabla f_{i,t} (x_{i,t}) ,x_{i, t}-y_{i,t} \rangle\\
	& + \langle\widetilde\nabla f_{i,t} (x_{i,t}) ,y_{i, t}-x_{i,t+1} \rangle\\
	&+\langle\widetilde\nabla f_{i,t} (x_{i,t})-\nabla f_{i,t} (x_{i,t}),x_{i, t}-x_t^{\star}    \rangle.
\end{align*}
Starting with the first term and via the boundedness of the $\widetilde\nabla f_{i,t} (x_{i,t})$, there holds that
\begin{align*}
    &\langle\widetilde\nabla f_{i,t}(x_{i,t}) ,x_{i, t+1}-x_t^{\star} \rangle=\frac{1}{\eta_t}\langle y_{i, t}-x_{i,t+1},x_{i, t+1}-x_t^{\star} \rangle \\
    &= \frac{1}{2\eta_t} (\|y_{i, t}-x_t^{\star}\|^2-\|y_{i, t}-x_{i,t+1}\|^2-\|x_{i, t+1}-x_t^{\star}\|^2)\\
    & \le \frac{1}{2\eta_t} (\|y_{i, t}-x_t^{\star}\|^2-\|x_{i, t+1}-x_t^{\star}\|^2) \\
    & \le \frac{1}{2\eta_t} (\|x_{i, t}-x_t^{\star}\|^2-\|x_{i, t+1}-x_{t}^{\star}\|^2)
\end{align*}
where in the last line we applied the convexity of the norm $\|\cdot\|$. 
The second term is bounded by the network error as
\begin{align*}
    &\langle\widetilde\nabla f_{i,t} (x_{i,t}) ,y_{i, t}-x_{i,t} \rangle  \\&=  \langle\widetilde\nabla f_{i,t} (x_{i,t}) ,y_{i, t}-\Bar{x}_t \rangle+ \langle\widetilde\nabla f_{i,t} (x_{i,t}) ,\Bar{x}_t-x_{i,t} \rangle\\
    &\le \lambda_t (\|x_{i,t}-\Bar{x}_t\|+
    \|y_{i, t}-\Bar{x}_t\|)\le 2\lambda_t \|x_{i,t}-\Bar{x}_t\|.
\end{align*}
The third term can be bounded by algorithm (\ref{algorithmonline}) 
\begin{align*}
    \langle\widetilde\nabla f_{i,t} (x_{i,t}) ,y_{i, t}-x_{i,t+1}\rangle=\eta_t\|\widetilde\nabla f_{i,t} (x_{i,t})\|^2\le \eta_t\lambda_t^2.
\end{align*}
Taking the summation and using Lemma \ref{lem_networkerr} lead to the validity of Lemma \ref{thm_1}.

\subsection{Proof of Lemma \ref{lem_highproboundonline}}\label{proofhighbound}
\textit{Proof of Lemma \ref{lem_highproboundonline}:} $\|\nabla f_{i,t}(x_{i,t})\|\le \frac{\lambda_t}{2}$ is already satisfied. First, we bound $\sum_{t=1}^{n} \eta_t\langle \theta_{i,t}^b,x_{i,t}-x^{*}\rangle$. By Assumption \ref{assu5} and Lemma \ref{lem_heavytaied}, it yields that
\begin{align*}
    &\sum_{t=1}^{T} \eta_t\langle \theta_{i,t}^b,x_{i,t}-x^{*}\rangle \le \sum_{t=1}^{T} \eta_t\|\theta_{i,t}^b\|\|x_{i,t}-x^{*}\|\\
    &\le 8B_X\sum_{t=1}^{T} \eta_t \lambda_t^{1-p} \sigma^p.
\end{align*}
We use the Freedman's inequality \cite{nguyen2023high} to bound another term. Note that $\{\langle \theta_{i,t}^u,x_{i,t}-x^{*}\rangle\}_{t\ge0}$ is a sequence of martingale difference. We denote
\begin{align*}
    &c_t=\eta_t\|\theta_{i,t}^u\|\|x_{i,t}-x^{*}\|=4B_X\eta_t\lambda_t,\\
    &\sum_{t=1}^{T} \mathbb{E} [(\eta_t\langle \theta_{i,t}^u,x_{i,t}-x^{*}\rangle)^2|\mathcal{F}_t]\\
    &\le 32\sigma^p B_X\sum_{t=1}^{T}\lambda_t^{2-p}\eta_t^2=F.
\end{align*}
Since $c_t\le c=8B_XB_g$ almost surely, it follows from the Freedman's inequality that
\begin{align*}
  &P\left( |\sum_{t=1}^{T} \eta_t\langle \theta_{i,t}^u,x_{i,t}-x^{*}\rangle|>a\text{ and } \right.\\
  & \left.\sum_{t=1}^{T} \mathbb{E} [(\eta_t\langle \theta_{i,t}^u,x_{i,t}-x^{*}\rangle)^2|\mathcal{F}_t]\le F \right)\\
  &\le 2\exp \left\{-\frac{a^{2}}{2F+\frac{2}{3}ca}\right\}.
\end{align*}
We choose $a$ such that
\begin{align*}
   2\exp \left\{-\frac{a^{2}}{2F+\frac{2}{3}ca}\right\} =\delta
\end{align*}
which gives
\begin{align*}
    a&=\frac{1}{2}\left(\frac{2c}{3}\log (\frac{2}{\delta})+\sqrt{\frac{4c^2}{9}\log^2 (\frac{2}{\delta})+8F}\right)\\
    &\le \frac{2}{3} c \log \frac{2}{\delta}+ \sqrt{2F}
\end{align*}
where the last inequality follows from the fact that $\sqrt{a+b}\le \sqrt{a}+\sqrt{b}$ for any $a,b\ge 0$.
\subsection{Proof of Theorem \ref{thm_dynamicregret}}\label{proofstatic}
\textit{Proof of Theorem \ref{thm_dynamicregret}:} Note that the exact gradient is bounded by $B_g$, one has that
\begin{align*}
    &f_t(x_{i,t})-f_t(x_t^{*})=f_t(x_{i,t})-f_t(\bar{x}_t)+f_t(\bar{x}_t)-f_t(x_t^{*})\\
    &\le \sum_{j=1}^{N} \langle \nabla f_{j,t} (x_{i,t}), x_{i,t}-\bar{x}_t\rangle+ f_t(\bar{x}_t)-f_t(x_t^{*})\\
    &\le B_gN\|x_{i,t}-\bar{x}_t\|+f_t(\bar{x}_t)-f_t(x_t^{*}).
\end{align*}
Moreover, we simplify the above as follows
\begin{align*}
    f_{j,t}(\bar{x}_t)-f_{j,t}(x_{i,t})    \le 
    B_g\|x_{i,t}-\bar{x}_t\|.
\end{align*}
Summing by $i\in [N]$, we derive the difference between the dynamic regret and the bound derived in Lemma \ref{thm_1}.
\begin{align*}
  &\sum_{i=1}^{N}f_t(x_{i,t})-f_t(x^{*}) \le \sum_{i=1}^{N} f_{i, t}\left(x_{i, t}\right)-f_{i, t}\left(x_t^{\star}\right )\\&+2NB_g\sum_{i=1}^{N}\|x_{i,t}-\bar{x}_t\|.
\end{align*}
With Assumption \ref{assu5}, we derive that
\begin{align*}
    \|x_{i, t}-x_t^{\star}\|^2-\|x_{i, t}-x_{t-1}^{\star}\|^2 \le 4B_X\|x_t^*-x_{t-1}^*\|
\end{align*}
where $B_X=\sup_{\mathcal{X}} \|x\|$. The inequality above indicates that
\begin{align*}
   & \frac{1}{2}(\|x_{i, t}-x_{t}^{\star}\|^2-\|x_{i, t+1}-x_{t}^{\star}\|^2)\\
   & \le \frac{1}{2} (\|x_{i, t}-x_{t-1}^{\star}\|^2-\|x_{i, t+1}-x_{t+1}^{\star}\|^2\\
   &+4B_X\|x_t^*-x_{t-1}^*\|).
\end{align*}
Thus, for any $i\in [N]$, we have that 
\begin{align*}
    &\sum_{t=1}^{T} \|x_{i, t}-x_{t}^{\star}\|^2-\|x_{i, t+1}-x_{t}^{\star}\|^2\\
    & \le \|x_{i, 1}-x_{1}^{\star}\|^2-\|x_{i, 2}-x_{2}^{\star}\|^2+ \sum_{t=2}^{T}(\|x_{i, t}-x_{t-1}^{\star}\|^2\\
    &-\|x_{i, t+1}-x_{t+1}^{\star}\|^2+4B_X\|x_t^*-x_{t-1}^*\|)\\
    &\le \|x_{i, 1}-x_{1}^{\star}\|^2+ 4B_XC_T.
\end{align*}
Rearranging terms in (\ref{thm_main_error}), we have for any $\delta\in (0,1)$, with the probability at least $1-\delta$, 
\begin{align*}
&\sum_{t=1}^{T}\eta_t(\sum_{i=1}^{N}f_{t}(x_{i,t})-f_t(x^{*})) \\&\le (5+2N\gamma\frac{1}{1-\beta})N\sum_{t=1}^{T}\eta_t^2\lambda_t^2\\
 &+(2B_X \lambda_1+2B_gN^2)N\gamma\eta_1\frac{1}{1-\beta}\\
 &+2B_gN^2(2+N\gamma\frac{1}{1-\beta})\sum_{t=1}^{T}\eta_t^2\lambda_t\\
 &+8N B_X \sigma^p \sum_{t=1}^{T}  \lambda_t^{1-p}\eta_t+2B_XC_T \\
 &+\frac{1}{2}\sum_{i=1}^{N}\|x_{i,1}-x_{1}^*\|^2+\frac{32}{3} B_X^2L \log \frac{2}{\delta}+ \sqrt{2F}
\end{align*}
where $F=32\sigma^p B_X\sum_{t=1}^{T}\lambda_t^{2-p}\eta_t^2$ and the inequality above is satisfied by the Fubini's theorem
\begin{align*}
    &\sum_{t=1}^{T}\lambda_t\eta_t\sum_{l=1}^{t-1}\beta^{t-l}\lambda_l\eta_l = \sum_{l=1}^{T-1}\sum_{t=l+1}^{T-1}\beta^{t-l}\lambda_t\eta_t\lambda_l\eta_l\\
    &=\sum_{l=1}^{T-1}\frac{1}{1-\beta} \lambda_{l+1}\eta_{l+1}\lambda_l\eta_l
    \le \sum_{l=1}^{T}\frac{1}{1-\beta} \lambda_l^2\eta_l^2.
\end{align*}
Due to the decreasing step size  $\eta_t$, the dynamic regret is bounded.
\subsection{Proof of Theorem \ref{thm_non-convex}}\label{proof_non-convex}
\textit{Proof of Theorem \ref{thm_non-convex}:} By the $L$-smoothness of $f_{i,t}$, there holds that
\begin{align*}
&f_{t}\left(\bar{x}_{t+1}\right) \le f_{t}\left(\bar{x}_{t}\right)+\frac{L}{2}\left\|\bar{x}_{t+1}-\bar{x}_{t}\right\|^{2}\\
&+\langle\nabla f_{t}\left(\bar{x}_{t}\right), \bar{x}_{t+1}-\bar{x}_{t}\rangle 
\\&=f_{t}\left(\bar{x}_{t}\right)+\frac{L}{2}\left\|\frac{\eta_{t}}{N} \sum_{i=1}^{N} \widetilde{\nabla} f_{i, t}\left(x_{i,t}\right)\right\|^{2}. 
\end{align*}
Considering $\langle\nabla f_{t}\left(\bar{x}_{t}\right), \frac{\eta _{t} }{N}\sum_{i=1}^{N} \widetilde{\nabla} f_{i, t}\left(x_{i, t}\right) \rangle$, we derive that
\begin{align*}
    &\langle\nabla f_{t}\left(\bar{x}_{t}\right), \frac{\eta _{t} }{N}\sum_{i=1}^{N} \widetilde{\nabla} f_{i, t}\left(x_{i, t}\right) \rangle-\frac{\eta_{t}}{N}\left\|\nabla f_{t}\left(\bar{x}_{t}\right)\right\|^{2}\\
    &=\frac{\eta _{t} }{N}\sum_{i=1}^{N}\Big(\langle\nabla f_{t}\left(\bar{x}_{t}\right), (\widetilde{\nabla} f_{i,t}\left(x_{i,t})-\nabla f_{i,t}\left(x_{i, t}\right)\right)\rangle\Big.
\\&\Big.-\langle\nabla f_{t}\left(\bar{x}_{t}\right),\left(\nabla f_{i,t}\left(x_{i, t}\right)-\nabla f_{i,t}\left(\bar{x}_{t}\right)\right)\rangle\Big).
\end{align*}
Note that $\left\|\widetilde{\nabla} f_{i, t}\left(x_{i, t}\right)\right\| \le \lambda_{t}$, it yields that
\begin{align*}
    &\eta_{t} \| \nabla f_{t}\left(\bar{x}_{t}\right) \|^{2} \le N\left(f_{t}\left(\bar{x}_{t}\right)-f_{t}\left(\bar{x}_{t+1}\right)\right)+\frac{L}{2} \eta_{t}^{2} \lambda_{t}^{2}\\
    &+\eta_{t} B_{g} L\left\|x_{i, t}-\bar{x}_{t}\right\|+\eta_{t} \sum_{i=1}^{N}\left\langle\theta_{i, t}, \nabla f_{t}\left(\bar{x}_{t}\right)\right\rangle .
\end{align*}
Next we provide the bound of the term $ \eta_{t} \left\langle\theta_{i, t}, \nabla f_{t}\left(\bar{x}_{t}\right)\right\rangle .$ Replacing
$2B_X$ to $Bg$ in Lemma \ref{lem_highproboundonline}, for any $\delta \in (0,1)$. with probability at least $1-\delta$
\begin{align*}
    &\sum_{t=1}^{T} \eta_{t} \sum_{i=1}^{N}\left\langle\theta_{i, t}, \nabla f_{t}\left(\bar{x}_{t}\right)\right\rangle\le 4N B_g \sum_{t=1}^{T} \eta_{t} \lambda t\left(\frac{\sigma}{\lambda_{t}}\right)^{p}\\
    &+\frac{8}{3} NB_{g}^{2} \log \frac{2}{\delta}+N\sqrt{2 F}
\end{align*}
where  $F=8 B_{g}^{2} \sigma^{p} \sum_{t=1}^{T} \lambda_{t}^{2-p} \eta_{t}^{2}$. In addition, note that
\begin{align*}
&\sum_{t=1}^{T}\left(f_{t}\left(\bar{x}_{t}\right)-f_{t}\left(\bar{x}_{t+1}\right)\right)-f_{1}\left(\bar{x}_{1}\right)+f_{T+1}\left(\bar{x}_{T+1}\right)\\
&=\sum_{t=1}^{T}\left(f_{t}\left(\bar x_{t}\right)-f_{t+1}\left(\bar x_{t+1}\right)\right)\le D_T .
\end{align*}
Note that for any $i\in [N]$ and $t\in[T]$, by $(a_1+\cdots+a_n)^2\le n(a^2+\cdots+a_n^2)$, one has that
\begin{align*}
    \|\nabla f_{t}(x_{i,t})\|^2 \le 2 \|\nabla f_{t}(\bar{x}_t)\|^2+2NL^2\|x_{i,t}-\bar{x}_t\|^2.
\end{align*}
Then by Lemma \ref{lem_networkerr}, Lemma \ref{cor_network^2}  and the decreasing property of $\eta _{t} $, the validity of the Theorem is ensured.

\end{document}